\documentclass[a4paper,12pt]{article}

\usepackage[margin=3cm,footskip=1cm]{geometry}

\usepackage[T1]{fontenc}
\usepackage{amsmath,amssymb,amsthm,bm,mathtools,xspace,booktabs,url}
\usepackage{graphicx,etoolbox}
\usepackage{color}

\usepackage{soul}
\usepackage[round]{natbib}

\usepackage[raggedright]{titlesec}

\numberwithin{equation}{section}

\theoremstyle{plain} 

\theoremstyle{definition} 

\usepackage{authblk}

\makeatletter
\newcommand\CorrespondingAuthor[1]{%
  \begingroup%
  \def\@makefnmark{}%
  \footnotetext{Corresponding author: #1}%
  \endgroup%
}
\makeatother

\makeatletter
\renewenvironment{abstract}{%
  \small%
  \providecommand\keywords{%
    \par\medskip\noindent\textit{Keywords:}\xspace}%
  \begin{center}%
    \bfseries \abstractname\vspace{-.5em}\vspace{\z@}%
  \end{center}%
  \quote%
}{\endquote}
\makeatother

\usepackage[all]{onlyamsmath}

\usepackage{lipsum}

\usepackage[utf8]{inputenc}
\usepackage{hyperref}
\usepackage{fixme}

\usepackage{subfigure}
\usepackage{float}

\newcommand{\gm}{\gamma}

\DeclareMathOperator\Var{Var}
\DeclareMathOperator\Cov{Cov}
\DeclareMathOperator\E{E}

\newcommand{\ta}{\theta}

\newcommand{\R}{\mathbb{R}}

\newcommand{\bx}{\mathbf{x}}
\newcommand{\by}{\mathbf{y}}

\newcommand{\bX}{\mathbf{X}}
\newcommand{\bT}{\mathbf{T}}
\newcommand{\bI}{\mathbf{I}}

\newcommand{\dd}{\,\mathrm{d}}

\begin{document}

\title{Towards optimal Takacs--Fiksel estimation}

\author[1]{Jean-François Coeurjolly}
\affil[1]{Laboratory Jean Kuntzmann, Univ. Grenoble Alpes,
  France, \texttt{Jean-Francois.Coeurjolly@univ-grenoble-alpes.fr}}

\author[2]{Yongtao Guan}
\affil[2]{Department of Management Science, University of Miami, USA, \texttt{yguan@bus.miami.edu}}

\author[3,4]{Mahdieh Khanmohammadi}
\affil[3]{Department of Computer Science, University of Copenhagen, Denmark}
\affil[4]{Department of Electrical Engineering and Computer Science,
  University of Stavanger, Norway, \texttt{mh.khanmohammadi@gmail.com}}

\author[5]{Rasmus Waagepetersen}
\affil[5]{Department of Mathematical Sciences, Aalborg University, Denmark, \texttt{rw@math.aau.dk}}

\date{\today}

\maketitle

\begin{abstract}
  The Takacs--Fiksel method is a general approach to estimate the
  parameters of a spatial Gibbs point process. This method embraces
  standard procedures such as the pseudolikelihood 
  and is defined via weight functions. In this
  paper we propose a general procedure to find weight functions which
  reduce the Godambe information and thus outperform pseudolikelihood
  in certain situations. The new procedure is applied to a
    standard dataset and to a recent neuroscience replicated point
    pattern dataset. Finally, the performance of the new
  procedure is investigated in a simulation study.

\end{abstract}

\section{Introduction}
Spatial Gibbs point processes are important models
for spatial dependence in point patterns \citep{lies00}
with a broad range of applications \citep[e.g.][]{stoyan:penttinen:00,illian:etal:08}.
Such processes are specified 
by a density with respect to a Poisson point process or, equivalently,
by the Papangelou conditional intensity. When the density or
Papangelou conditional intensity has a parametric form, popular
options for parameter estimation include maximum likelihood
\cite[e.g.][]{ogattane84,penttinen:84,geyer:99,moeller:waagepetersen:04}, maximum
pseudolikelihood
\cite[e.g.][]{besag:77,jensen:moeller:91,goulard:sarkka:grabarnik:96,baddeley:turner:00,billiot:coeurjolly:drouilhet:08},
maximum logistic regression likelihood \citep{baddeley:etal:14} and
Takacs-Fiksel estimation \cite[e.g.][]{fiksel:84,takacs:86,tomppo:86,coeurjolly:etal:12}.

Maximum likelihood estimation for a Gibbs point process
requires computationally intensive estimation of an unknown
normalizing constant in the density function. This explains why
alternative estimation methods have been studied. Takacs-Fiksel
estimation is an estimating function method based on the general
Georgii-Nguyen-Zessin integral equation involving the Papangelou
conditional intensity and a user-specified weight function. A particular
choice of the weight function recovers the score of the
pseudolikelihood.
Pseudolikelihood estimation has an intuitively appealing
motivation and is by far the most popular estimation method in
practical applications of Gibbs point processes with a user-friendly implementation in the \texttt{spatstat}
package \citep{baddeley:turner:05}. Logistic regression likelihood
estimation for Gibbs point processes was recently introduced to eliminate
a bias problem coming from the Berman-Turner approximation of
the pseudolikelihood
\citep{berman:turner:92,baddeley:turner:00}. The logistic
  regression can be viewed as a  computationally efficient
approximation of the pseudolikelihood. Hence in the following, we may  not
differentiate between the pseudolikelihood and the logistic regression methods.

There are infinitely many weight functions that can be used to obtain
Takacs-Fiksel estimates and usually  weight functions are chosen by
ad hoc reasoning paying attention to ease of implementation or to handle patterns which are sampled in situ, like in \citet{tomppo:86}.
In the small scale simulation study in \citet{coeurjolly:etal:12},
pseudolikelihood and Takacs-Fiksel methods are compared for the
  Strauss model. The weight functions used for the Takacs-Fiksel method
  resulted in explicit expressions for the parameter estimates
    (no optimization required). In this simulation study, smallest
    standard errors were obtained with the pseudolikelihood method. Based on a
larger simulation study, \cite{diggle:etal:94}
concluded that pseudolikelihood provided reasonable and robust
estimates for any point process model and boundary condition
considered although significant bias was observed in cases of point processes models with strong
interaction. For the Takacs-Fiksel method, \cite{diggle:etal:94} used  weight functions leading to an
  estimation procedure similar to minimum contrast estimation using the empirical $K$-function, see e.g. \cite{waagepetersen:guan:09}.
The authors noted that this type of Takacs-Fiksel
estimation gave poor results relative to pseudolikelihood, in
particular for point process models with weak
  interaction. \cite{diggle:etal:94} concluded
  that the bias was small for both Takacs-Fiksel and pseudolikelihood
  in case of small-to-medium strength of interaction. Based on the existing literature it is
not clear whether pseudolikelihood is an optimal Takacs-Fiksel method
in terms of minimizing estimation variance and it is
  interesting to investigate whether and when the Takacs-Fiksel method can outperform pseudolikelihood. 

In this paper our aim is to develop a systematic approach to construct a weight function that can lead to  more efficient estimation than existing methods. Our approach
is motivated by the one considered in
\cite{guan:jalilian:waagepetersen:15} who considered estimation of the
intensity function of a spatial point process and identified the
optimal estimating function within a class of
estimating functions based on the Campbell formula
\citep[e.g.][]{moeller:waagepetersen:04}. Their optimal estimating
function was derived from a sufficient condition equating the sensitivity matrix for the optimal estimating function
and the covariance between the optimal estimating function and an
arbitrary estimating function.

Extending the ideas in \cite{guan:jalilian:waagepetersen:15} to Gibbs
 point processes is not straightforward. One problem is that
covariances of Takacs-Fiksel estimating functions are not available
in closed forms (see Section~\ref{sub:comparison} for more details). 
For this reason our new weight
function only approximately satisfies the aforementioned sufficient
condition. Nevertheless, we show in a simulation study that the new
weight function may yield better estimation accuracy and is closer to
fulfilling the sufficient condition than pseudolikelihood. Another
issue is that the practical implementation of the new method is more
computationally demanding than the pseudolikelihood, especially for point
patterns of high cardinality. 

The rest of the paper is organized as
follows. Section~\ref{sec:methodology} gives  background on Gibbs
point processes and Takacs-Fiksel estimation and presents our new
methodology. In Section~\ref{sec:examples} we apply the methodology to the Spanish towns dataset as well as
a recent replicated point pattern dataset from neuroscience. 
Motivated by these two examples, Section~\ref{sec:sim} presents a
simulation study. Details of implementation are given in Appendix~\ref{sec:implementation}.

\section{Background and methodology} \label{sec:methodology}

\subsection{Gibbs point processes}\label{sec:Gibbs}

A point process $\bX$ on $\R^d$ is a locally finite random subset
of $\R^d$, meaning that $\bX \cap B$ is finite for every bounded $B \subset \R^d$. In this
paper we assume that $\bX$ is confined to and observed on a bounded
region $W \subset \R^d$ so that $\bX$ becomes a finite point process taking values in $\Omega$, the set of finite point configurations in $W$.

The distribution of $\bX$ is assumed to be specified by a parametric
density $f(\cdot ;\theta):\Omega \to [0,\infty)$ with respect to the
Poisson process of unit intensity. The density is of the form
\begin{equation}\label{eq:gibbsdensity}
  f(\by;\theta) \propto H(\by) e^{V(\by;\theta)}
\end{equation}
where $\theta \in \Theta \subseteq \R^p$ is a $p$-dimensional
parameter vector, $H:\Omega \to [0,\infty)$ serves as a baseline or
reference factor, and $V:\Omega \to\R$ is often called the
potential. For the Strauss model, for example, $H(\by)=1$, $\theta=(\theta_1,\theta_2)$ and  \[ V(\by;\theta) = \theta_1 n(\by) + \theta_2 s_R(\by) \]
where $n(\by)$ is the cardinality of $\by$ and for $R>0$, $s_R(\by)$ is the number of
subsets $\{u,v\}$, $u,v \in \by$, of pairs of $R$-close neighbours in
$\by$. Thus values of $\theta_2<0$ promote point configurations with
few $R$-close neighbours. The
Strauss hard core model is the modification of the Strauss model where
$H(\by)$ is one if all interpoint distances are greater than some hard core
distance $0 \le \delta<R$ and zero otherwise. Thus for the Strauss hard
core model, all points must be separated by a distance greater than~$\delta$.

Assuming that $H$ is hereditary, i.e. $H(\by \cup u)>0$ implies $H(\by)>0$ for any $u\in W$ and $\by \in \Omega$, the Papangelou conditional intensity of $\bX$ exists and is defined by
\[
  \lambda (u,\by;\theta) = \frac{f(\by \cup u; \theta)}{f(\by;\theta)} = H(u,\by)e^{V(u,\by;\theta)}
\]
where $H(u,\by)=\mathbf 1 \{H(\by)>0 \}H(\by \cup u)/H(\by)$ and
$V(u,\by;\theta)=V(\by\cup u;\theta)-V(\by;\theta)$. For the Strauss model,
$H(u,\by)=1$ and $V(u,\by;\theta)=\theta_1+\theta_2 s_R(u,\by)$ where
$s_R(u,\by)$ is the number of points in $\by$ with distance to $u$
  less  than $R$. For the Strauss hard core model, $H(u,\by)$ is one if all the points of
$\by\cup u$ are separated by distances greater than $\delta$ and zero
otherwise. In this paper we assume finite range, i.e.\
for some $0<R<\infty$ and for any $u\in W$ and $\by\in \Omega$
\begin{equation}
  \label{eq:FR}
  \lambda(u,\by; \theta) = \lambda\{u , \by\cap B(u,R);\theta\}
\end{equation}
where $B(u,R)$ is the ball with center $u$ and radius
$R$. Thus the conditional intensity of a point $u$ given
$\by$ only depends on the $R$-close neighbours in $\by$. This is
obviously satisfied for the Strauss and Strauss hard core models. Intuitively, $\lambda(u, \bx;\theta)
\dd u$  is the conditional probability that a point of $\bX$ occurs in
a small neighbourhood $B_u$ of volume $\dd u$ around the location $u$, given
$\bX$ outside $B_u$ is equal to $\bx$; see \cite{georgii:76} for a
general presentation and \cite{coeurjolly:moeller:waagepetersen:15}
for links with Palm distributions.

Note that
$\lambda$ and $f$ are in one-to-one correspondence. Hence the
distribution of $\bX$ can equivalently be specified in terms of the
Papangelou conditional intensity. Gibbs point processes can also be
characterized  through the Georgii-Nguyen-Zessin formula
\citep[see][]{georgii:76,nguyen:zessin:79}, which states that for any $h:W \times \Omega\to \R$ (such that the following expectations are finite)
 \begin{equation}
   \label{eq:GNZ} \E \sum_{u \in \bX} h(u,\bX\setminus u) = \E \int_{W} h(u,\bX) \lambda(u,\bX;\theta) \dd u.
 \end{equation}

Conditions ensuring the existence of Gibbs point processes constitute
a full research topic \citep[see e.g.][and the references
therein]{dereudre:drouilhet:georgii:12}. In case of a finite Gibbs
point process existence is equivalent to that the right hand side of
\eqref{eq:gibbsdensity} can be normalized to be a probability density.
We here just note that the
Strauss model exists whenever
$\theta_1 \in \R$ and $\theta_2 \le 0$ while the Strauss hard core
model exists for all $\theta_1,\theta_2\in \R$ when $\delta>0$. Other examples of
point process models can be found in \cite{moeller:waagepetersen:04}; see also Section~\ref{sec:vesicles} for an example of an inhomogeneous model.



\subsection{Takacs-Fiksel estimation}

Let $h=(h_1,\ldots,h_q)^\top$ where $h_i:W \times \Omega\to \R$, $i=1,\ldots,q$ are real functions parameterized by
$\theta$  where $q$ is greater than or equal to the dimension $p$ of
$\theta$. Takacs-Fiksel estimation \citep{fiksel:84,takacs:86} is based on the
Georgii-Nguyen-Zessin formula~\eqref{eq:GNZ} which implies that 
\begin{equation}\label{eq:hi} \E \sum_{u\in \bX} h_i(u,\bX\setminus u;\theta) - \E\int_W
h_i(u,\bX;\theta)\lambda(u,\bX;\theta)\dd u =0, \quad i=1,\ldots,q. 
\end{equation}
In the original formulation of Takacs-Fiksel estimation, an estimate
of $\ta$ was obtained by minimizing the sum of squares
\[
\sum_{i=1}^q \left\{ \sum_{u\in \bX} h_i(u,\bX\setminus u;\theta) -\int_W h_i(u,\bX;\theta)\lambda(u,\bX;\theta)\dd u
  \right\}^2.
\]
Considering the case $p=q$, a related approach is to define a $p$-dimensional estimating function $e_h$ by
\begin{equation}\label{eq:ef}
e_h(\theta)  = \sum_{u\in \bX} h(u,\bX\setminus u;\theta) -\int_W h(u,\bX;\theta)\lambda(u,\bX;\theta)\dd u.
\end{equation}
In the point process literature, the term Takacs-Fiksel estimate is also
used for an estimate obtained by solving the estimating equation
$e_h(\ta)=0$ with respect to $\ta$ \citep[see e.g.\ ][]{lies00,baddeley:rubak;turner:15}). This is the type of Takacs-Fiksel
estimation considered in this paper.
By \eqref{eq:hi} the function $e_h(\ta)$ is unbiased, i.e., $\E e_h(\theta)=0$.
The
  asymptotic properties of the Takacs-Fiksel estimate for a general weight
  function $h$ have been established by~\citet{coeurjolly:etal:12}
  including a derivation of the asymptotic covariance matrix.

With the weight function $h(u,\by;\theta)=\dd \log \lambda(u,\by;\theta) / \dd
\ta$, \eqref{eq:ef} is the score function for the log pseudolikelihood function. The corresponding estimate
can be obtained using standard statistical software and its
statistical properties have been deeply studied in the literature
\citep[e.g.\ ][]{jensen:moeller:91,mase97,jenskuen94,billiot:coeurjolly:drouilhet:08,baddeley:etal:14}. However, the
pseudolikelihood score has
not been shown to be optimal within the class of Takacs-Fiksel
estimating functions.  In
  the following section our aim is to construct a competitor to the pseudolikelihood in terms of
statistical efficiency.

\subsection{Towards optimality}

We begin this section by reviewing key quantities related to
  estimating functions. We refer to \citet{heyde:97} or
  \citet[Section~2.3]{guan:jalilian:waagepetersen:15} for further details.
For an estimating function $e_h$, the sensitivity matrix $S_h$ is
defined as $S_h= -\E \left\{
  \frac{\dd}{\dd \theta^\top} e_h(\theta)\right\}$ and the
  covariance matrix of the estimating function is $\Sigma_h=\Var\{e_h(\theta)\}$. From these
the Godambe information matrix is obtained as
\[ G_h=S_h^\top \Sigma_h^{-1} S_h. \]
In applications of estimating functions, the inverse Godambe matrix
provides the approximate covariance matrix of the associated parameter
estimate. An estimating function $e_\phi$ is said to be Godambe optimal in a class of
estimating functions $e_h$ indexed by a set $C$ of functions $h$, if the
difference $G_\phi-G_h$ is non-negative definite for all $h \in
C$. Following \cite{guan:jalilian:waagepetersen:15}, a sufficient
condition for $e_\phi$ to be optimal is that for every estimating
function $e_h$, $h \in C$,
\begin{equation}\label{eq:suffcond} \Cov\{e_h(\theta),e_\phi(\theta)\}=S_h.\end{equation}

In the context of Takacs-Fiksel estimating functions \eqref{eq:ef},
\begin{align*}
  \frac \dd{\dd \theta^\top} e_h(\theta)=& \sum_{u\in \bX} \frac{\dd}{\dd \theta^\top}h(u,\bX\setminus u;\theta)-\int_{W} \left\{\frac{\dd}{\dd \theta^\top} h(u,\bX;\theta) \right\}\lambda(u,\bX;\theta) \dd u \\
  &-\int_W h(u,\bX;\theta) \frac{\dd}{\dd \theta^\top } \lambda(u,\bX;\theta) \dd u.
\end{align*}
So by the Georgii-Nguyen-Zessin formula~\eqref{eq:GNZ},
\begin{equation}
  \label{eq:Sf}
  S_h= \E \int_W h(u,\bX;\theta) \frac{\dd}{\dd \theta^\top } \lambda(u,\bX;\theta) \dd u.
\end{equation}
The definition of the estimating function $e_h$ actually corresponds
to the concept of innovations for spatial point processes
\citep{baddeley:etal:05}. \citet{coeurjolly:rubak:13} investigated the
problem of estimating the covariance between two innovations which,
here, corresponds to the covariance between two estimating
functions. Assuming the right hand side below is finite,
\citet[Lemma~3.1]{coeurjolly:rubak:13} established
that
\begin{align}
\Cov\{&e_h(\theta) , e_g(\theta) \} =\E \Big [ \int_{W} h(u,\bX;\theta)g(u,\bX;\theta)^\top\lambda(u,\bX;\theta) \dd u \nonumber\\
&+ \int_W \int_W h(u,\bX;\theta) {g(v,\bX;\theta)^\top} \left\{ \lambda(u,\bX;\theta)\lambda(v,\bX;\theta)-\lambda(\{u,v\},\bX;\theta)\right\} \dd u \dd v \nonumber\\
&+  \int_W\int_W\Delta_v h(u,\bX;\theta)\Delta_u g(v,\bX;\theta)^\top
\lambda(\{u,v\},\bX;\theta)\dd u \dd v \Big ] \label{eq:covefeg}
\end{align}
where for any $u,v\in W$ and any $\by\in \Omega$, the second order Papangelou conditional intensity $\lambda(\{u,v\},\by;\theta)$ and the difference operator $\Delta_u h(v,\by;\theta)$ are given by
\begin{align*}
\lambda(\{u,v\},\by;\theta) &= \lambda(u,\by)\lambda(v, \by \cup u;\theta) =\lambda(v,\by;\theta)\lambda(u, \by \cup v;\theta) \\
\Delta_u h(v,\by;\theta) &= h(v ,\by \cup u;\theta) - h(v,\by;\theta).
\end{align*}

Returning to the condition \eqref{eq:suffcond}, we introduce for any $\by \in \Omega$ the operator $T_\by$ acting on $\R^p$ valued functions $g$,
\begin{equation} \label{eq:operator}
  T_\by g(u) = \int_{W} g(v) t(u,v,\by;\theta) \dd v,
\end{equation}
where
\begin{equation}\label{eq:kernel}
  t(u,v,\by;\theta) = \lambda(v,\by;\theta) \left\{ 1-\frac{\lambda(v,\by \cup u;\theta)}{\lambda(v,\by;\theta)}\right\}.
\end{equation}
The finite range property of the Papangelou conditional intensity
implies that for any $v\notin B(u,R)$, $t(u,v,\by;\theta)=0$. So the domain
of integration in \eqref{eq:operator} is actually just $W \cap
B(u,R)$. From~\eqref{eq:Sf}-\eqref{eq:covefeg}, \eqref{eq:suffcond} is
equivalent to $\E(A)+ \E(B)=0$ 
where
  \begin{align}
  A&=  \int_W h(u,\bX;\theta) \lambda(u,\bX;\theta) \left\{
    \phi(u,\bX;\theta) - \frac{\lambda^{(1)}(u,\bX;\ta)}{\lambda(u,\bX;\theta)} + T_{\bX}\phi (u,\bX;\theta)\right\}^\top \dd u\label{eq:A} \\
B&=  \int_W \int_{W} \Delta_v h(u,\bX;\theta)\Delta_u \phi(v,\bX;\theta)^\top \lambda(\{u,v\},\bX;\theta)\dd u \dd v \label{eq:B}
  \end{align}
and $\lambda^{(1)}(u,\bX;\ta)=\dd \lambda(u,\bX;\theta)/\dd \ta$.

The expectation $\E (B)$ is very difficult to evaluate. Moreover, in
the context of asymptotic covariance matrix estimation for the
pseudolikelihood, \citet{coeurjolly:rubak:13} remarked that the
contribution of the term~\eqref{eq:B} to the covariance
$\Cov\{e_h(\theta),e_\phi(\theta)\}$ was negligible. In the
following we will neglect the term~\eqref{eq:B} and call
`semi-optimal' a function $\phi:W\times \Omega\to \R^p$ (parameterized by $\theta$) such that
for any $h:W\times \Omega \to \R^p$, $\E(A)=0$. We discuss this choice
in more detail in Section~\ref{sub:comparison}. Considering
\eqref{eq:A}, $\phi$ is semi-optimal if for
any  $\by\in \Omega$, $\phi(\cdot,\by;\theta)$ is the solution to the
  Fredholm integral equation \citep[e.g.\ chapter~3 in][]{hackbusch:95}
 \begin{equation}
   \label{eq:phiOpt}
      \phi(\cdot,\by;\theta)  +
      T_\by \phi(\cdot,\by;\theta) = \frac{\lambda^{(1)}(\cdot,\by;\theta)}{\lambda(\cdot,\by;\theta)}.
 \end{equation}
In practice, this equation is solved numerically, see
Section~\ref{sec:outlineimplementation} and
Appendix~\ref{sec:implementation}.

Having solved \eqref{eq:phiOpt}, 
the
covariance  and sensitivity matrices for the resulting estimating
function
  \begin{equation}\label{eq:optimalef}
  e_\phi(\theta) =\sum_{u\in \bX}\phi(u,\bX\setminus u;\theta)-\int_W
  \phi(u,\bX;\theta)\lambda(u,\bX;\theta)\dd u
\end{equation}
are given by
\begin{align*}
  S &= \E \int_{W} \phi(u,\bX;\theta) \lambda^{(1)}(u,\bX;\theta)^\top \dd u\\
  \Sigma &= S+ \E \int_W \int_{W} \Delta_u \phi(v,\bX;\theta) \Delta_v \phi(u,\bX;\theta)^\top \lambda(\{u,v\},\bX;\theta) \dd u \dd v.
\end{align*}
Note that for a truly optimal $\phi$, we would have $S=\Sigma$. In the
simulation studies in Section~\ref{sec:sim_strauss}, we
  investigate for the Strauss model how close $S$ and $\Sigma$ are for our semi-optimal $\phi$.

\subsection{A comparison with optimal intensity estimation} \label{sub:comparison}

In this section we compare our approach with the problem of optimal intensity estimation for spatial point processes investigated and solved by \cite{guan:jalilian:waagepetersen:15}. They consider estimating equations of the form
\begin{equation}\label{eq:ehguanetal}
  e_h(\theta) = \sum_{u\in \bX} h(u;\theta) -\int_W h(u;\theta) \lambda(u;\theta) \dd u  
\end{equation}
to estimate a parametric model $\lambda(\cdot;\theta)$ for the intensity
function. Using the Campbell theorem, the sufficient
condition~\eqref{eq:suffcond} for $e_\phi$ to be optimal is equivalent to
\begin{equation}\label{eq:guanetal}
  \int_W h(u;\theta) \lambda(u;\theta) \left\{ 
\phi(u;\theta) - \frac{\lambda^{(1)}(u;\theta)}{\lambda(u;\theta)} + T \phi(u;\theta)
  \right\} \dd u  =0,
\end{equation}
where the operator $T$ acting on $\R^p$ valued functions $g$ is given by
\[
  T g(u) = \int_W g(v) \lambda(v;\theta) \{ \text{pcf}(v-u)-1 \} \dd v,
\]
where $\text{pcf}$ is the so-called pair correlation function. We now compare~\eqref{eq:guanetal} with the condition $\E(A)+\E(B)=0$ where $A$ and $B$ are given by~\eqref{eq:A}-\eqref{eq:B}.

The terms $A$ and~\eqref{eq:guanetal} are  similar where the roles played by the first and second order
Papangelou conditional intensities  in $A$ correspond to the ones played by
the intensity and pair correlation functions
in~\eqref{eq:guanetal}. However, solving $\E (A)=0$ is much more complex
than solving~\eqref{eq:guanetal}. This is mainly due to that
the function $\phi$ is, in the present paper, a function of both $u \in
W$ and $\by \in \Omega$: the estimating function~\eqref{eq:ehguanetal}
requires the computation of $\phi(v)$ for $v \in W$ while our
estimating function~\eqref{eq:ef} requires the evaluation of
$\phi(u,\bx\setminus u)$ for data points $u \in \bx$ as well as the
evaluation of $\phi(v,\bx)$ for $v \in W$, see
Appendix~\ref{sec:implementation} for details on the numerical implementation.

Further, due to the simpler form of~\eqref{eq:ehguanetal}, there is no
term like $B$ appearing in~\eqref{eq:guanetal}. In the problem we
consider, it is not possible to include $B$ and still obtain a
condition that is manageable in practice. As mentioned in the previous
section, $B$ has been observed to be negligible in certain
applications in which case our approach should lead to an estimating equation
able to outperform the pseudolikelihood method. As far as we know, our
work is the first in this direction.

\subsection{Outline of implementation}\label{sec:outlineimplementation}

Let $\bx=\{x_1,\ldots,x_n\}$ denote an observation of $\bX$.
In practice we approximate the integral in \eqref{eq:operator} by
numerical quadrature whereby the integral equation \eqref{eq:phiOpt}
becomes a matrix equation where for each $\by=\bx, \bx\setminus
x_1,\dots,\bx \setminus x_n$ the unknown quantity is the vector with
 components $\phi(u_j;\by;\theta)$ where $u_j$, $j=1,\ldots,m,$ denote
 the $m$ quadrature points. The detailed description of the
 implementation which is quite technical and heavy on notation is
 given in the Appendix. 

 Depending on the choice of $m$ the solution of
 the matrix equation involving a sparse $m \times m$ matrix can be
 computationally quite costly. Since this has to be done for each $\by=\bx, \bx\setminus
x_1,\dots,\bx \setminus x_n$, our method is not recommended for point
patterns of large cardinality. In the data examples we provide some
examples of computing times. In our applications we use quadrature
points located on a regular grid covering $W$. Thus e.g.\ for a $50
\times 50$ grid, $m=2500$.

The solution of the matrix equation
depends  on that the matrix involved is positive definite. For purely
repulsive models like the Strauss model and the Strauss hard core
model with $\ta_2<0$ this was always the case in the data examples and
the simulation studies. However in cases of positive interaction where
it is possible for the integral operator kernel \eqref{eq:kernel} to
be negative, the condition of positive definiteness was sometimes
violated. This e.g.\ happened for the multiscale process considered in
Section~\ref{sec:vesicles} and Section~\ref{sec:simulation_inh}. In
such cases we simple returned the pseudolikelihood estimate or
restricted the kernel  \eqref{eq:kernel} to be non-negative, see
Section~\ref{sec:vesicles} and Section~\ref{sec:simulation_inh} for
more details.

\section{Data examples} \label{sec:examples}

In this section we present two applications. For the first Spanish towns dataset, we apply a Strauss hard core model. For the second replicated point pattern dataset from neuroscience, we consider a multiscale hard core model.

\subsection{Spanish towns  dataset}

\begin{figure}[htbp]
\centering
\includegraphics[scale=.45]{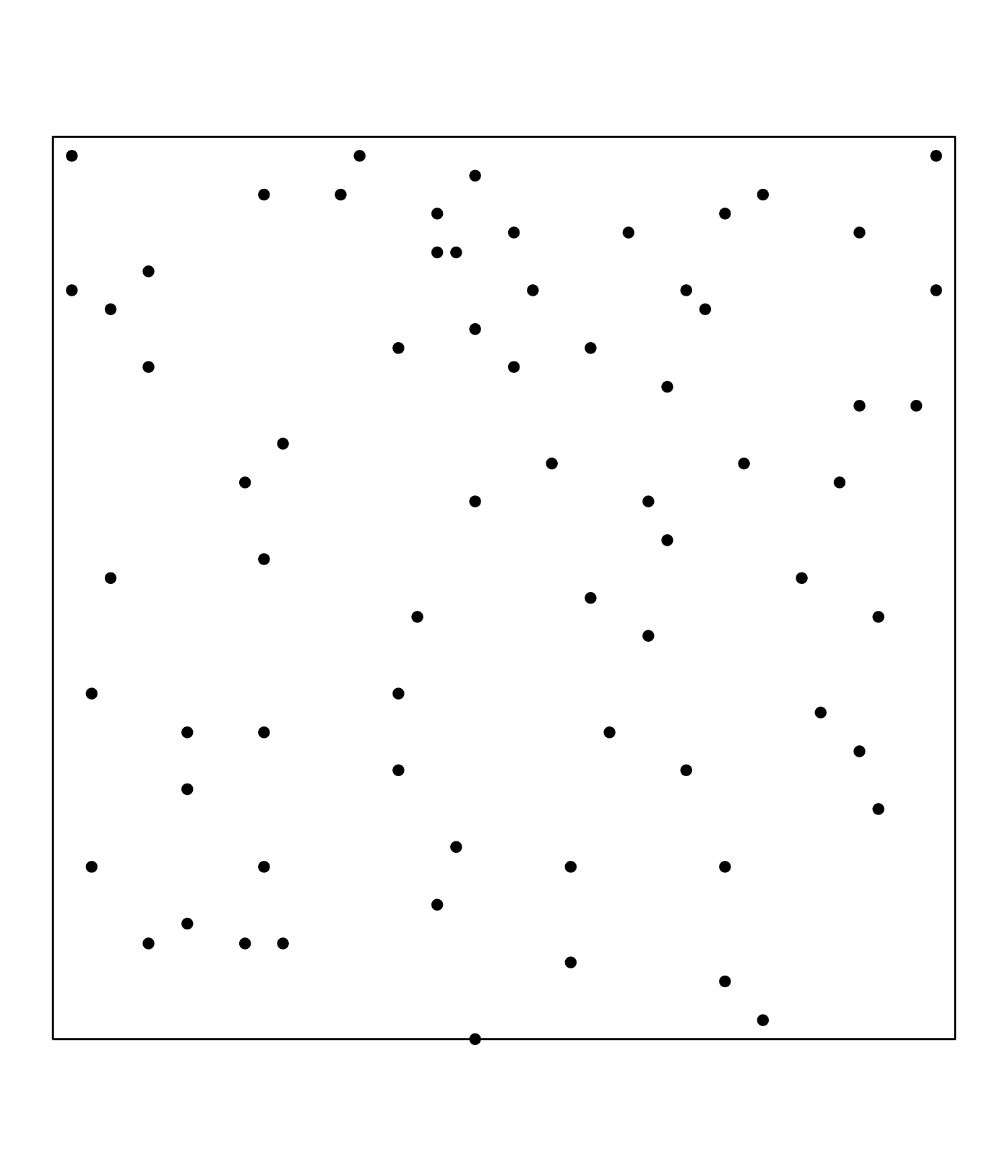}
\includegraphics[scale=.5]{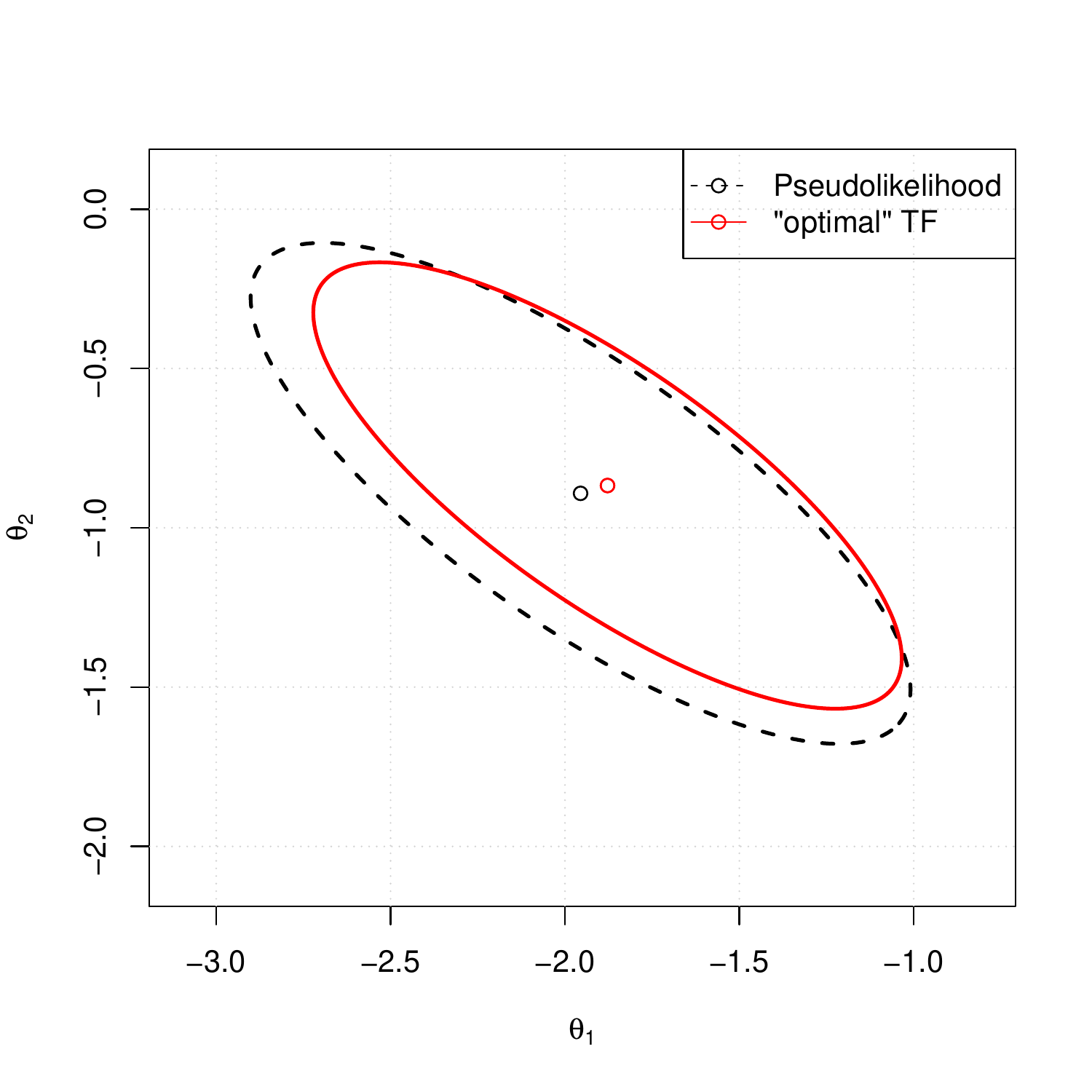}
\caption{\label{fig:towns} Left: locations of 69 Spanish towns in a 40
miles by 40 miles region. Right: 95\% confidence ellipses for the parameters $(\theta_1,\theta_2)$ for the pseudolikelihood and the semi-optimal Takacs-Fiksel methods.}
\end{figure}

\cite{ripley:88} and subsequently
\cite{illian:etal:08} proposed to model the Spanish towns
dataset (see Figure~\ref{fig:towns}) using a Strauss hard core
model.
We {compare results regarding estimation of $\theta_1$ and
$\theta_2$ using respectively the new semi-optimal method and the
pseudolikelihood. For the hard core distance and the interaction range
we use the values $\widehat \delta=0.83$ and $\widehat R=3.5$ obtained
using maximum likelihood by~\citet[p.~170]{illian:etal:08}.
For the pseudolikelihood we use the unbiased
logistic likelihood implementation introduced in
\cite{baddeley:etal:14} with a stratified quadrature point process
 on a $50\times50$ grid. The same grid is used for the semi-optimal
Takacs-Fiksel method.

Pseudolikelihood and semi-optimal estimates are presented in
Table~\ref{tab:towns}. The computing time for obtaining the semi-optimal estimates is 25 seconds on a Lenovo W541 laptop.
\begin{table}[ht]
\begin{center}
\begin{tabular}{lcc}
\hline
& $\theta_1$ & $\theta_2$\\
\hline
SO & -1.88 (0.12) & -0.87 (0.08) \\
PL &  -1.96 (0.15) & -0.89 (0.10) \\
Ratio se  & 0.79& 0.79 \\
\hline
\end{tabular}
\caption{\label{tab:towns} Results for Strauss hard core model applied to the
  Spanish towns dataset. First and second row: semi-optimal (SO)
  and pseudolikelihood (PL) estimates with estimated standard errors
  in parantheses. Last row: ratio of semi-optimal standard errors to
  pseudolikelihood standard errors. }
\end{center}
\end{table}
Standard errors of the  pseudolikelihood estimates (respectively the semi-optimal Takacs-Fiksel estimates) are estimated  from 500 simulations of the  model fitted using the pseudolikelihood (respectively the semi-optimal procedure). The standard errors are clearly reduced with the semi-optimal Takacs-Fiksel method. In
addition to these results, we compute the Frobenius norm of the
estimated covariance matrices for the pseudolikelihood and optimal
Takacs-Fiksel estimates (which actually correspond to
estimates of the inverse Godambe matrices). We obtained the
norm values 0.25 and
0.20 for the pseudolikelihood and the semi-optimal Takacs-Fiksel
methods respectively.
Based on the asymptotic normality results
established by \cite{baddeley:etal:14} and \cite{coeurjolly:etal:12},
we construct 95\% confidence ellipses for the parameter vector
$(\theta_1,\theta_2)$. The ellipses are depicted in the right plot of
Figure~\ref{fig:towns}. The area of the confidence region for the
semi-optimal Takacs-Fiksel method is  81\% of the one
for the pseudolikelihood. In line with the simulation study in
  Section~\ref{sec:sim_strauss},
these empirical findings show that more precise parameter estimates
can be obtained with the semi-optimal Takacs-Fiksel method compared to pseudolikelihood.

\subsection{Application to synaptic vesicles} \label{sec:vesicles}

Synapses are regions in the brain where nerve impulses are
transmitted or received. Inside the synapses, neurotransmitters are
carried by small membrane-bound compartments called synaptic
vesicles. Recently \cite{khanmohammadi:etal:14} studied whether stress
affects the spatial distribution of vesicles within the synapse. The
data used for the study originated from microscopial images of slices
of brains from respectively a group of 6 control rats and a group of 6
stressed rats. The images were annotated to identify
the boundaries of the synapses and possible mitochondria in the
synapses, the locations of the vesicles, and the extents of the so-called active zones where the
vesicles release their contents of
neurotransmitters. Figure~\ref{fig:annotation} shows annotations of
two images. One synapse was
  annotated for each rat except for one control rat where two synapses
  were annotated. Thus in total 7 synapses from the control rats and 6 synapses
from the stressed rats were annotated. For each
synapse several images corresponding to several slices of the synapse
were annotated. We restrict here attention to the middle slice for
each synapse. Thus, our data consist of 7 and 6 annotated images for
respectively the control and the stressed rats. The side lengths of
enclosing rectangles for the synapses range between 396 and 1009
  nm with a mean of 663 nm. Further details
on the dataset can be found in \cite{khanmohammadi:etal:14}.

\begin{figure}[htbp]
\centering
\begin{tabular}{ll}
\includegraphics[scale=.37]{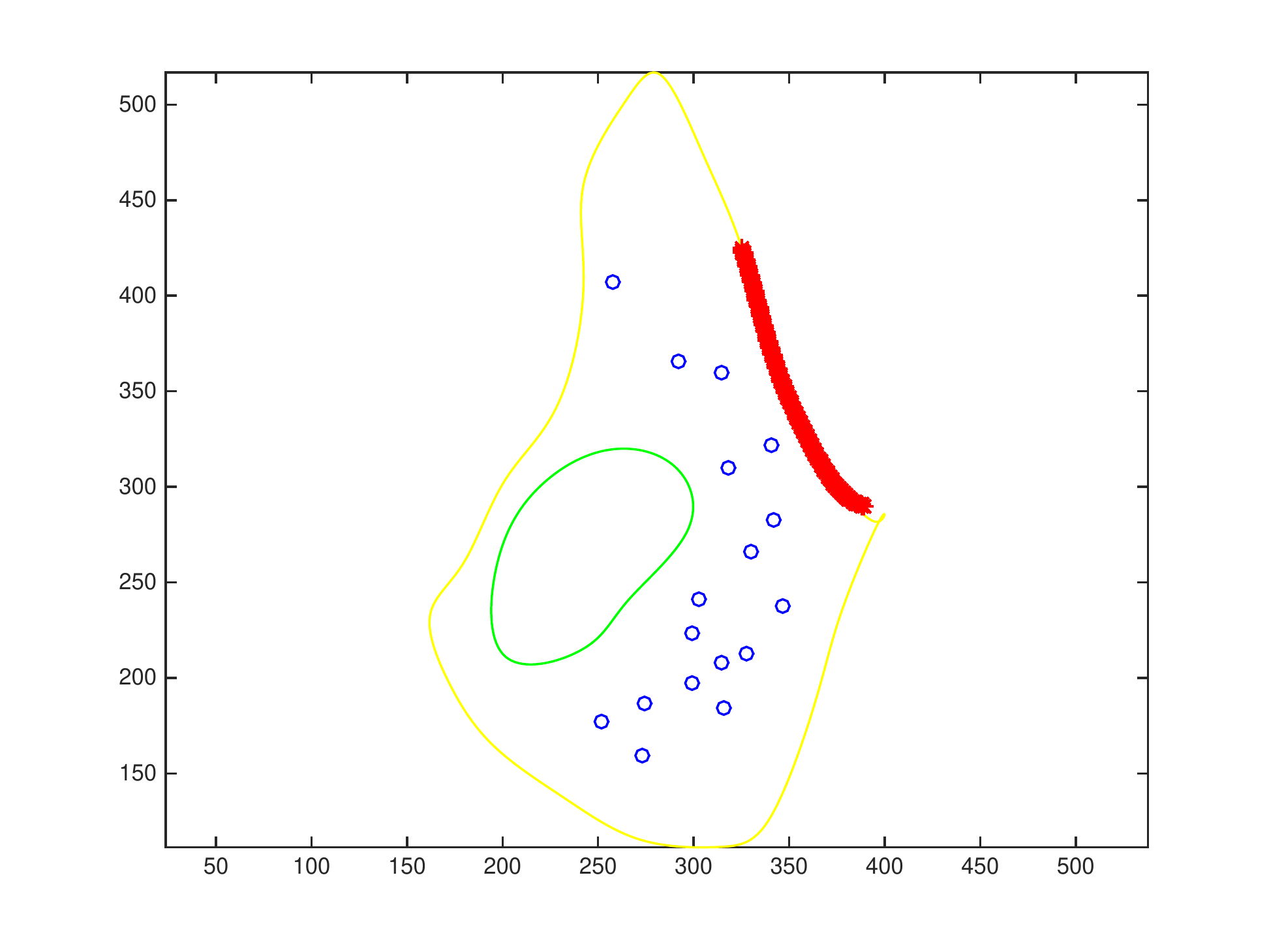} &
\includegraphics[scale=.37]{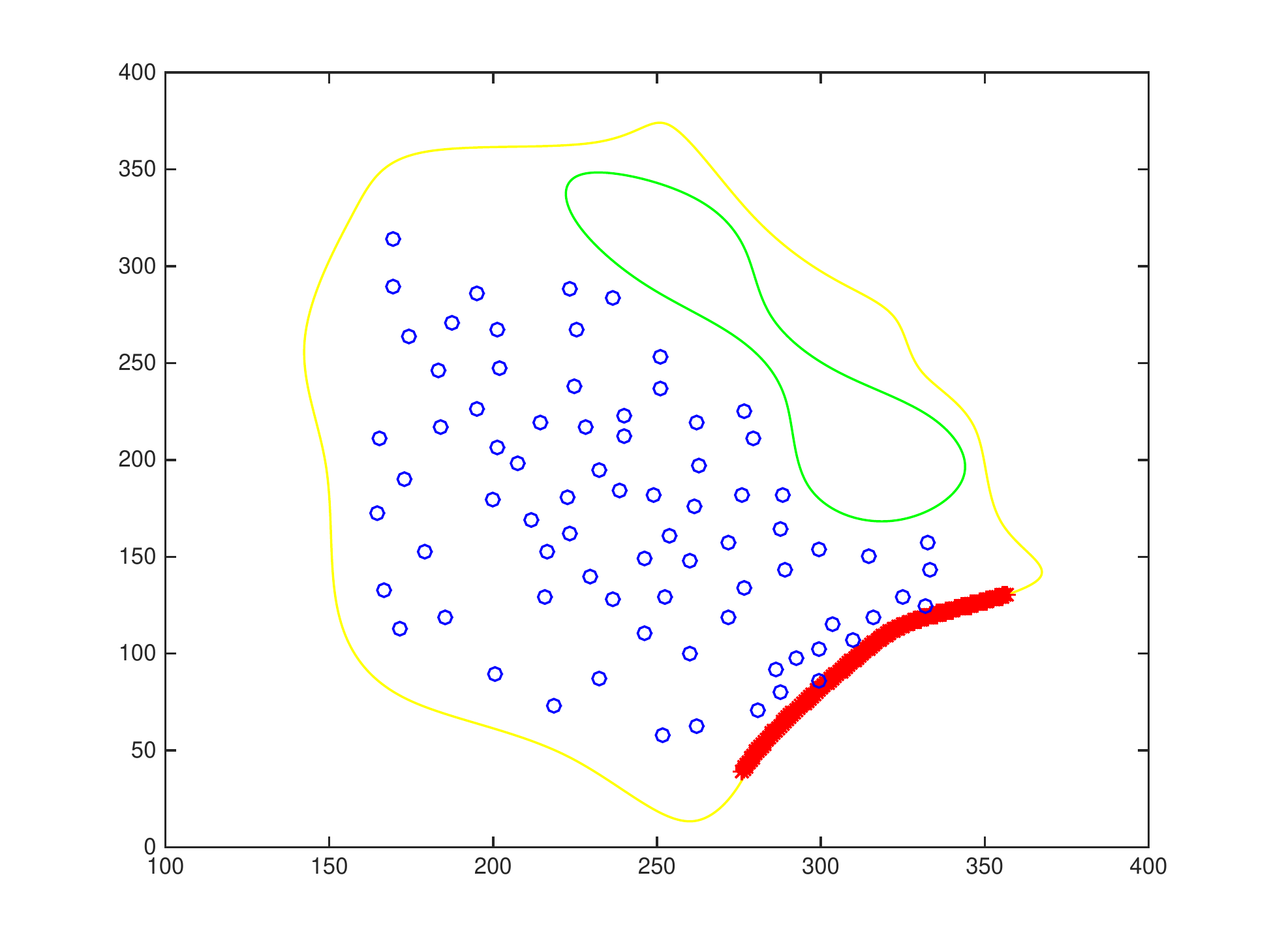}
\end{tabular}
\caption{\label{fig:annotation} Plots of vesicle locations for a control (left) and a stressed (right) synapse. The active zones are shown by the red curves while the
  green curves show the boundaries of mitochondria. The yellow curves
  show the boundaries of the synapses given by their cell membranes.}
\end{figure}

\subsubsection{Point process models for locations of vesicles}

For the $i$th image of type $t=C,S$ (control or stressed) we
consider the centers of the vesicles as a realization of a finite spatial
point process $\bX_{ti}$ with observation window $W_{ti}$ defined by
the boundary of the synapse excluding areas
occupied by mitochondria. We further assume that the pairs $(W_{ti},\bX_{ti})$
of the observation windows and the spatial point
processes for different images are independent, and that pairs $(W_{ti},\bX_{ti})$ of the same type $t$ are
identically distributed. \cite{khanmohammadi:etal:14} modelled the locations of the vesicles as an
inhomogeneous Strauss hard core process, but noted some evidence of
aggregation of vesicles at a larger scale not accommodated by this
model. After considering the Strauss hard core model for reference, we therefore
extend it to a
multiscale model with an additional interaction
term.
More precisely, for a
location $u$ and a configuration $\bx$ of vesicle
locations in a synapse of type $t$, the conditional
intensity is of the form
\begin{equation}
 \lambda_t(u,\bx;\theta) = \exp\{\theta_{0t} +\theta_{1t} d(u)+
 \theta_{2t} s_r(u,\bx)+ \theta_{3t}
 s_{r,R}(u,\bx)\}H_\delta(u,\bx),\quad t=C,S
\label{eq:conddens}
\end{equation}
where $0 \le \delta< r< R$, $\theta_{0t},\theta_{1t},\theta_{2t},\theta_{3t} \in \R$, $d(u)$ is the distance to the
active zone, $s_r(u,\bx)$ is the number of points in $\bx$ with
distance to $u$ smaller than $r$, $s_{r,R}(u,\bx)$ denotes the number
of points in $\bx$ with distance to $u$ in the interval $[r,R]$, and the hard core term $H_\delta(u,\bx)$ is
one if all points in $\bx\cup u$ are separated by a distance greater than $\delta$ and
zero otherwise. The inhomogeneous Strauss hard  core model is the special case with
$\theta_{3t}=0$. As in \cite{khanmohammadi:etal:14}, the hard core distance $h$ is set to 17.5nm
corresponding to the average diameter of a vesicle and the
interaction distance $r$ is set to 32.5nm. We further choose
the value of $R=107.5$ nm to maximize a profile
pseudolikelihood based on the data for all 13 images. We finally
scale the distances $d(\cdot)$ by a factor $10^{-3}$ in order to obtain
$\theta_{lt}$ estimates of the same order of magnitude.

\subsubsection{Inference for replicated point patterns}\label{sec:pooled}

Let $e^{(ti)}$ denote an estimating function (either the
pseudolikelihood score or the semi-optimal) for the $ti$th image and denote by $\theta_t$ the
vector of parameters to be inferred, $t=C,S$. Then for each type
$t=C,S$ we optimally form a pooled estimating function $e^{(t)}$ as the
sum of the $e^{{(ti)}}$'s. Following standard asymptotic arguments for pooled
independent estimating functions, the variance of the
corresponding parameter estimate is approximated as
\[ \Var {\big\{}\sqrt{n_t}(\hat \theta_t -\theta_t) {\big\}} \approx (S^{t1})^{-1} \Var {\big\{ e^{(t1)} \big\}} (S^{t1})^{-1} \]
\cite[the so-called sandwich estimator,][]{song:07} where $n_t$ is the number of replicates of type $t$ and $S^{t1}$ is the sensitivy matrix associated with $e^{{(t1)}}$. In practice we replace $S^{t1}$ and
${\Var\{e^{(t1)} \}} $ by their empirical estimates replacing the unknown $\theta_t$
with its estimate. We conduct the pseudolikelihood estimation
for the replicated data using the user-friendly \texttt{mppm}
procedure in the \texttt{R} package \texttt{spatstat} while using our
own code to evaluate the approximate variances of the
pseudolikelihood or semi-optimal estimates.

\subsubsection{Results for synapse data}

Table~\ref{tbl:strausshc} shows parameter estimates and associated
standard errors obtained for the Strauss hard core model with either
the pseudolikelihood or the
semi-optimal approach.
\begin{table}
\centering
\begin{tabular}{l|ccc}
\hline
     & $\theta_{0C}$ & $\theta_{1C}$ & $\theta_{2C}$ \\
SO &   -6.74 (0.21) &  -0.68 (0.20) &  -2.79 (1.00) \\
PL &   -6.80 (0.21) &  -0.66 (0.22) & -2.61 (1.00) \\
Ratio se & 1.01 & 0.80& 1.02 \\ \hline
& $\theta_{0S}$ &     $\theta_{1S}$ & $\theta_{2S}$ \\
SO &  -5.51 (0.23) & -1.14 (0.47) &  -0.61 (0.21) \\
PL&   -5.63 (0.24) & -1.12 (0.49) &  -0.45 (0.24)\\
Ratio se &  0.93& 0.91 & 0.78 \\
\hline
\end{tabular}
\caption{Results for Strauss hard core process. Semi-optimal (SO) and pseudolikelihood (PL) estimates with estimated standard
errors in parantheses and ratios of semi-optimal standard
errors to pseudolikelihood standard errors. First three rows:
control. Last three rows: stressed. Ratios of standard errors were
computed before rounding standard errors to two digits.}\label{tbl:strausshc}
\end{table}
Except for $\theta_{2S}$, the semi-optimal and pseudolikelihood
  estimates are fairly similar. The qualitative conclusions based on
the two types of estimates are identical: negative dependence of conditional intensity
on distance and repulsion between vesicles both in the control and the stressed
group. The estimated semi-optimal standard
errors are smallest for all parameters in the stressed group and
$\theta_{1C}$ in the control groups. For $\theta_{0C}$ and
$\theta_{2C}$ the estimated semi-optimal and pseudolikelihood standard
errors are very similar.

Due to the aforementioned evidence of large
scale aggregation we turn to the multiscale model for a more detailed
comparison of the control and the stressed group. {To avoid problems with
negative definiteness in the numerical implementation (see
Section~\ref{sec:outlineimplementation}) we restricted
$\theta_{2t}$ and $\theta_{3t}$, $t=C,S,$ to be less or equal to zero inside the
kernel \eqref{eq:kernel}}. Parameter estimates
and associated standard errors for the multiscale model are shown in
Table~\ref{tbl:multiscale}. 
\begin{table}[htbp]
\centering
\begin{tabular}{l|cccc}
\hline
  & $\theta_{0C}$ & $\theta_{1C}$ & $\theta_{2C}$ & $\theta_{3C}$ \\
SO &   -8.87  (0.12) &  -0.47  (0.11) &  -3.40  (0.78) &  0.45 (0.02) \\
PL & -8.71  (0.13)& -0.46  (0.12)& -3.35  (0.80)&  0.40  (0.02)  \\
Ratio se &  0.94 & 0.93 & 0.98 & 1.01  \\ \hline
& $\theta_{0S}$ & $\theta_{1S}$ & $\theta_{2S}$ & $\theta_{3S}$ \\
SO & -7.87  (0.21) & -0.48  (0.26) &  -1.17  (0.18)&  0.26  (0.02)\\
PL & -7.84  (0.19) & -0.69  (0.27) & -1.12  (0.18) &  0.24  (0.02) \\
Ratio se & 1.12 & 0.96 & 0.97 &1.07\\
\hline
\end{tabular}
\caption{\label{tbl:multiscale} Results for multiscale
  process. Semi-optimal (SO) and pseudolikelihood (PL) estimates with estimated standard
errors in parantheses and ratios of semi-optimal standard
errors to pseudolikelihood standard errors. First three rows:
control. Last three rows: stressed. Ratios of standard errors were
computed before rounding standard errors to two digits.}
\end{table}
The qualitative conclusions based on the pseudolikelihood and the
semi-optimal estimates coincide.
All parameters are significantly
different from zero (assuming estimate divided by its
standard error is approximately $N(0,1)$). In particular the positive
estimates of $\theta_{3C}$ and $\theta_{3S}$ confirm that there is
aggregation at a larger scale.
The negative estimates of $\theta_{1t}$, $t=C,S$ indicate that the conditional
intensity is decreasing as a function of distance to the active zone
while, according to the estimates of $\theta_{2t}$, there appears to be a strong respectively moderate small scale interaction
for the control group and the stressed group.

To test the hypotheses: $H_l$:
$\theta_{lS}=\theta_{lC}$, $l=0,1,2,3$, we consider statistics of the form $(\hat
\theta_{lS}-\hat \theta_{lC})/\sqrt{\text{se}_{lS}^2+\text{se}_{lC}^2}$
where  $\text{se}_{lt}^2$ denotes the estimated standard error of
the estimate $\hat \theta_{lt}$ $t=S,C$ and the estimates are
obtained either using the semi-optimal approach or pseudolikelihood. Under the hypothesis this
statistic is approximately $N(0,1)$. According to these tests, $H_0$
and $H_2$ and $H_3$ are rejected while $H_1$ is not irrespective of the
estimation method ($p$-values adjusted for multiple testing using \cite{holm:79}'s procedure are $<0.001$, 0.96, 0.02, $<0.0001$
and $<0.001$, $0.88$, $0.02$, $<0.0001$ 
for semi-optimal and pseudolikelihood, respectively). There is thus evidence that the
repulsion between vesicles is stronger for the control
rats than for the stressed rats which means that the vesicles tend to
form more regular patterns for the control rats.

For the control data,
the estimated standard errors are smallest with the semi-optimal
approach except for $\theta_{3C}$. For the stressed rats the
semi-optimal standard errors are smallest for $\theta_{1S}$ and
$\theta_{2S}$ but not for $\theta_{0S}$ and $\theta_{3S}$. Overall, a
clear pattern is not visible. Note also that the ratios of estimated
standard errors should be interpreted with care as they are obviously
subject to sampling error. We further compare semi-optimal and
pseudolikelihood for the multiscale model in Section~\ref{sec:simulation_inh}.

The total computing time for fitting the multiscale model to the control and stressed rats datasets is respectively 24 and 18 minutes on a Lenovo W541 laptop. The results shown are obtained using a $80 \times 80$ grid. For comparison we also tried a $60 \times 60$ grid and got very similar results but with computing time reduced to 6 and 5 minutes. The computing time can be further reduced by running computations for each synapse in parallel.

\section{Simulation study} \label{sec:sim}

This section investigates by simulation studies the  performance
of the semi-optimal approach when applied to the models used in the
data examples in Section~\ref{sec:examples}.

\subsection{Strauss model} \label{sec:sim_strauss}

{To study the performance of our semi-optimal estimating function
relative to the pseudolikelihood score, we apply both estimating
functions to simulations of a Strauss process
(Section~\ref{sec:Gibbs}) on the unit square (equivalent to a
  Strauss hard core model with $\delta=0$)}. We  use the
\texttt{spatstat} \citep{baddeley:turner:05} procedure \texttt{rStrauss()} to generate exact
simulations of the Strauss process for  $\beta=\exp(\ta_1)=100$ and all combinations of
$R=0.04,0.08,0.12$ and $\gamma=\exp(\ta_2)=0.1,0.2,0.4,0.8$. In the
following, we specify parameter settings in terms of the parameters
$(\beta,\gamma)$ due to their ease of interpretation. We however
compare the efficiency of the pseudolikelihood and the semi-optimal
estimation methods based on the estimates of $\ta_1$ and $\ta_2$ which are
typically closer to being normal than the estimates of $\beta$ and
$\gamma$. 

The semi-optimal estimating function is implemented using a
$50\times50$ or a $75\times 75$ grid. As in the previous sections, for the pseudolikelihood we use the unbiased logistic likelihood implementation with a stratified quadrature
  point process. For both estimation methods, $R$ is assumed to be known
 and equal to the value used to generate the simulations. The biases of the semi-optimal and the
    pseudolikelihood estimates are very similar. The bias relative to
    the true parameter estimate is very small ($-1$ to $1$\%) for
    $\theta_1$ but can be substantial for $\theta_2$ (up to 28\% in
    case of $\gamma=\exp(\theta_2)=0.8$ and $R=0.04$). However, in all
    cases the bias is
    negligible relative to the estimation variance. We therefore
    focus on root mean square error (RMSE) when comparing the two
    methods. Due to the small bias relative to estimation variance we would have obtained essentially the same results by
replacing RMSE with the standard error of the simulated
estimates.

For the parameters $\theta_1$
and $\theta_2$, Table~\ref{tbl:rmse} shows for each
parameter setting, the RMSE
of the pseudolikelihood estimates minus the RMSE of the semi-optimal estimates relative to the
RMSE of the pseudolikelihood estimates. For each parameter
  setting, the RMSEs are estimated from a sample of 1000 parameter estimates obtained
  from 1000 simulations. We omit simulations where all interpoint distances are
larger than $R$ and thus neither the pseudolikelihood nor the
semi-optimal estimates of $(\theta_1,\theta_2)$ exist (if one considered instead $(\beta,\gamma)$, an alternative would be
to let the estimate of $\gamma$ be equal to zero when all interpoint
distances are greater than $R$). The  reported relative differences in
RMSE are subject to Monte
Carlo error. Table~\ref{tbl:rmse} therefore also shows estimated standard
errors for these obtained by applying a bootstrap
to each Monte Carlo sample.
\begin{table}[H]
\centering
\begin{tabular}{r|ll|ll|ll|ll}
 \hline
& \multicolumn{8}{c}{Interaction parameter $\gamma=\exp(\theta_2)$} \\
Range   & \multicolumn{2}{c}{0.1 }& \multicolumn{2}{c}{0.2 }& \multicolumn{2}{c}{0.4 }& \multicolumn{2}{c}{0.8} \\
and grid & $\ta_1$ &$\ta_2$ &$\ta_1$ &$\ta_2$ &$\ta_1$ &$\ta_2$ &$\ta_1$ &$\ta_2$ \\
\hline
\multicolumn{1}{r|}{\small $R=0.04$} &&&&&&&\\
(50,50) & -1 (0.5) & 1 (0.3) & -1 (0.4) & 0 (0.2) & -2 (0.3) & 0 (0.3) & -2 (0.2) & 0 (0.3) \\
  (75,75) & -1 (0.4) & 0 (0.2) & -1 (0.3) & 0 (0.2) & -1 (0.3) & 0 (0.2) & -1 (0.2) & 0 (0.2) \\
  \multicolumn{1}{r|}{\small$R=0.08$}&&&&&&& \\
 (50,50) & 3 (1) & 1 (0.6) & 0 (0.9) & 0 (0.6) & 5 (0.8) & 2 (0.7) & 4 (0.5) & 2 (0.7) \\
   (75,75) & 2 (0.8) & 2 (0.5) & 3 (0.9) & 2 (0.6) & 2 (0.8) & 1 (0.7) & 1 (0.5) & 1 (0.6) \\
     \multicolumn{1}{r|}{\small$R=0.12$} &&&&&&&\\
         (50,50) & 3 (1.2) & 5 (0.8) & 4 (1.1) & 3 (1) & 6 (1.3) & 4 (1.3) & 5 (1) & 4 (1.2) \\
    (75,75) & 3 (1.3) & 5 (0.9) & 7 (1.2) & 4 (1) & 7 (1.3) & 5 (1.2) & 7 (1.1) & 6 (1.4) \\
 \hline
   \end{tabular}\caption{RMSE for pseudolikelihood minus RMSE for semi-optimal relative to RMSE
for pseudolikelihood (in percent) in case of estimation of
$\theta_1$ and $\theta_2$.  Grids of $50\times 50$ or $75\times 75$ quadrature points are considered. Numbers between brackets are bootstrap standard errors (in percent).
 } \label{tbl:rmse}
\end{table}

In case of $R=0.04$, there is no
efficiency improvement by  using the semi-optimal estimating function. In fact the semi-optimal approach appears to be
sometimes slightly worse than pseudolikelihood, even when taking into
account the Monte Carlo error of the estimated relative differences in RMSE.
However, for $R=0.08$ and $R=0.12$ the semi-optimal approach is
  always better with decreases up to
5 and 7\% in RMSE for semi-optimal relative to RMSE for
pseudolikelihood. The results regarding the performance of
  semi-optimal relative to pseudolikelihood are quite similar for the
  two choices of grids with slightly more favorable results for
  semi-optimal in case of the $75 \times 75$ grid and $R=0.12$. For $R=0.04$ the expected number of points varies from 71 to 92
  when $\gamma$ varies from 0.1 to 0.8. For $R=0.08$ and $R=0.12$ the
  corresponding numbers are 43-75 and 28-61. There does not seem to be
  a clear dependence on the
  expected number of simulated points regarding the performance of
  semi-optimal relative to pseudolikelihood.

For an optimal estimating function and at the true parameter value,
the covariance matrix $\Sigma$ of the estimating function coincides with the
sensitivity matrix $S$.  Table~\ref{tbl:optimal}
shows the estimated relative differences $(S_{ij}-\Sigma_{ij})/\Sigma_{ij}$ in
percent for the pseudolikelihood and semi-optimal estimating
functions.
\begin{table}[h]
\centering
\begin{tabular}{r|rrr|rrr|rrr|rrr}
\hline
& \multicolumn{12}{c}{Interaction parameter $\gamma=\exp(\theta_2)$} \\
Interaction range & \multicolumn{3}{c}{0.1 }& \multicolumn{3}{c}{0.2 }& \multicolumn{3}{c}{0.4 }& \multicolumn{3}{c}{0.8} \\
\hline
 SO, $R=0.04$& 3 & -17 & -48 & 5 & -5 & -49 & 2 & 15 & -44 & 3 & -4 & -46 \\
 0.08 & -5 & -31 & -50 & -2 & -16 & -44 & 3 & -7 & -40 & 1 & 1 & -33  \\
  012& -10 & -38 & -45 & -10 & -36 & -49 & -7 & -23 & -41 & -5 & -10 & -26 \\
 \hline
 PL, $R=0.04$ & 41 & -4 & -46 & 40 & 11 & -47 & 27 & 34 & -41 & 12 & 4 & -43 \\
 0.08 & 87 & -4 & -46 & 80 & 20 & -35 & 71 & 35 & -27 & 30 & 27 & -20  \\
0.12  & 119 & -3 & -37 & 102 & 6 & -35 & 87 & 31 & -18 & 41 & 32 & 0 \\
\hline
\end{tabular}\caption{Relative difference (in percent) between sensitivity and
  variance of estimating function, $(S_{ij}-\Sigma_{ij})/\Sigma_{ij}$,
  $ij=11,12,22$. Upper three rows: semi-optimal (SO). Lower three rows: pseudolikelihood (PL). A grid of $75\times 75$ quadrature points is used.}\label{tbl:optimal}
\end{table}
For neither of the estimating functions, the covariance and
sensitivity agree. However, in general the relative deviations are
larger for pseudolikelihood than for semi-optimal. To investigate
  this further, we have computed the Frobenius norms of the matrices with entries  $(S_{ij}-\Sigma_{ij})/\Sigma_{ij}$. The results (not shown) indicate that the Frobenius norms for the semi-optimal method are smaller than the corresponding ones for the pseudolikelihood for all cases except $\gamma=0.8$ and $R=0.04$.


{
\subsection{Multiscale hard core model} \label{sec:simulation_inh}

For the multiscale hard core process with conditional intensity
\[
 \lambda(u,\bx;\theta) = \exp\{\theta_{0} +\theta_{1} d(u)+ \theta_{2} s_r(u,\bx)+ \theta_{3} s_{r,R}(u,\bx)\}H_\delta(u,\bx)
\]
we consider three settings inspired by the synaptic vesicles data example. The
process is simulated on the unit square and we let $d(u)=x$ be the first coordinate of $u=(x,y) \in [0,1]^2$ and
$\theta_1=-0.5$ to obtain a decreasing trend in the first spatial
coordinate. The terms $s_r$ and $s_{r,R}$ and the factor $H_\delta$
are defined as below \eqref{eq:conddens} with hard core distance
$\delta=0.01$ and interaction ranges $r=0.08$ and $R=0.16$. The
values of $\gm_2=\exp(\theta_2)$ and $\gm_3=\exp(\theta_3)$ are
$(\gm_2,\gm_3)=(0.2,1)$, $(\gm_2,\gm_3)=(0.2,1.25)$ and
$(\gm_2,\gm_3)=(0.2,1.5)$. Note that $\gm_3=1$ corresponds to the null hypothesis of no large scale interaction. The values of
$\beta=\exp(\theta_0)=100,60,40$ are adjusted to produce on average 40
points for each of the three parameter settings. 

Following the
approach in the previous section, Table~\ref{tbl:rmsemultiscale} shows for each
parameter setting, the root mean square error (RMSE)
of the pseudolikelihood estimates minus the RMSE of the semi-optimal estimates relative to the
RMSE of the pseudolikelihood estimates (estimated using 1000
simulations and using a $60 \times 60$ grid). As in the previous
section we omit cases where neither the semi-optimal nor the
pseudolikelihood estimates of $\theta=(\ta_0,\ta_1,\ta_2,\ta_3)$
exist. For the setting with $\gm_3=1.5$ (last two rows of
Table~\ref{tbl:rmsemultiscale}) we often (24\% of the simulations)
encountered negative
definiteness in the numerical computations. The results in the last
row of Table~\ref{tbl:rmsemultiscale} are obtained with
$\theta_2$ and $\theta_3$ restricted to be less than or equal to
zero inside the kernel \eqref{eq:kernel} in which case the
negative definiteness only occurred for 2\% of the
simulations. As mentioned in Section~\ref{sec:outlineimplementation},
we just return the pseudolikelihood estimates in the cases where the
semi-optimal approach fails.}
\begin{table}[!htb]
\centering
   

   \begin{tabular}{lllll}
\hline
& $\ta_0$ & $\ta_1$ & $\ta_2$ & $\ta_3$\\
\multicolumn{1}{l}{$(\beta,\gm_2,\gm_3)$} \\
\hline
(100,0.2,100) \\
&5 (1) & 6 (1) & 3 (1) & 10 (1)\\
(60,0.2,1.25) \\ 
& 7 (1) & 4 (1) & 2 (1) & 10 (1)  \\ \hline
(40,0.2,1.5) \\
&3 (1) &1 (1) & -1 (1) & 1 (1)\\
(40,0.2,1.5) \\ 
& 0 (1) &5 (1) &2 (1) &5 (1)\\
\hline
\end{tabular}
\caption{RMSE for pseudolikelihood minus RMSE for semi-optimal relative to RMSE
for pseudolikelihood (in percent) in case of estimation of
$(\theta_0,\theta_1,\theta_2,\theta_3)$.  Numbers between brackets are
bootstrap standard errors (in percent). In the last row,
$\theta_2$ and $\theta_3$ are restricted to be non-positive in the
kernel \eqref{eq:kernel}.
 } \label{tbl:rmsemultiscale}
\end{table}

For $\gm_3=1$ and $\gm_3=1.25$ the semi-optimal approach
  performs much better than pseudolikelihood with decreases in relative
  RMSE of up to 10\%. In case of $\gm_3=1.5$ decreases up to 5\% are
  obtained in the last row with the restricted version of the kernel
  \eqref{eq:kernel}. With $\gm_3=1.5$ and without restriction (third
  row), only very moderate decreases are obtained for
  $\ta_0,\ta_1$ and $\ta_3$ while for $\ta_2$, the semi-optimal approach appears to be slightly worse than pseudolikelihood.

\section{Discussion}

In this paper, we have investigated the scope for outperforming the
pseudolikelihood by tuning weight functions for the Takacs-Fiksel
estimator. Due to the complicated nature of moments for Gibbs point
processes, the method is less straightforward than the one proposed
by~\cite{guan:jalilian:waagepetersen:15} for estimating the intensity
function of a spatial point process. Therefore our new Takacs-Fiksel
method is not guaranteed to be optimal. It is also computationally
more expensive than the approach in
\cite{guan:jalilian:waagepetersen:15} because the weight function needs
to be evaluated for both the observed point pattern and all patterns
obtained by omitting one point at a time.

In the simulation studies we have demonstrated that the new semi-optimal approach can
yield better statistical efficiency both for purely repulsive point
processes and also certain multiscale models with moderate positive
interaction. In presence of stronger positive
interactions, the method is, however, susceptible to  numerical
problems. These numerical problems can, to some extent, be mitigated by
introducing restrictions in the kernel \eqref{eq:kernel}. When comparing the methods, the higher computational complexity of the semi-optimal approach
  should also be taken into account. Thus, while we have made a significant
  step towards optimal Takacs-Fiksel estimation there is still room
  for further improvement.

\subsubsection*{Acknowledgments}
We thank the editor and two referees for detailed and
constructive comments that helped to improve the paper. R. Waagepetersen is supported by the Danish Council for Independent Research | Natural
Sciences, grant 12-124675, "Mathematical and Statistical Analysis of
Spatial Data", and by the "Centre for Stochastic Geometry and Advanced
Bioimaging", funded by grant 8721 from the Villum Foundation. \text{J.-F.} Coeurjolly is supported by ANR-11-LABX-0025 PERSYVAL-Lab (2011, project OculoNimbus).

\bibliographystyle{royal}
\bibliography{opt}

\appendix 

\section{Implementation} \label{sec:implementation}

Let $\bx=\{x_1,\dots,x_n\}$ denote a realization of $\bX$.
To solve $e_\phi(\theta)=0$, we use Newton-Raphson iterations starting
at the pseudolikelihood estimate with the
Hessian matrix estimated by the empirical sensitivity matrix $\hat S =  \int_{W}
\phi(u,\bx;\theta) \lambda^{(1)}(u,\bx;\theta)^\top \dd u$. To evaluate
$e_\phi$ and $\hat S$ we need to solve \eqref{eq:phiOpt} with respect
to $\phi(\cdot;\by;\theta)$ for all
$\by=\bx, \bx\setminus x_1,\dots,\bx \setminus x_n$.

\subsection{Symmetrization}

To ease the implementation and in particular the use of Cholesky
decompositions, we symmetrize the operator $T_\by$. This is possible if
we assume for any $u,v\in W$ and $\by\in \Omega$, the ratio
$\lambda(v,\by\cup u;\theta)/\lambda(v,\by;\theta)$ is symmetric in
$u$ and $v$. This assumption  is valid for instance for all pairwise interaction point processes. Indeed, the Papangelou conditional intensity of such processes is given by
$\lambda(u,\by ;\theta) = e^{\sum_{w\in \by} \psi(\{w,u\};\ta)}$ where
$\psi$ is a real valued function, whereby
$\lambda(v,\by\cup u;\theta)/\lambda(v,\by;\theta) =
e^{\psi(\{v,u\};\ta)}$.

We now multiply each term of \eqref{eq:phiOpt} by
$\sqrt{\lambda(\cdot,\by;\theta)}$ and reformulate the problem to
solve
\begin{equation}
  \label{eq:phiOpt2}
  \widetilde \phi(\cdot,\by;\theta)   + \widetilde T_\by \widetilde \phi (\cdot,\by;\theta) =
   \frac{\lambda^{(1)}(\cdot,\by;\theta)}{\sqrt{\lambda(\cdot,\by;\theta)}}
\end{equation}
with respect to the function $\widetilde
\phi(\cdot,\by;\theta)=\sqrt{\lambda(\cdot,\by;\theta)}\phi(\cdot,\by;\theta)$ where
$\widetilde T_\by$ is the operator with kernel
\[  \widetilde t(u,v,\by;\theta)= \sqrt{\lambda(u,\by;\theta)\lambda(v,\by;\theta)}\left\{ 1-\frac{\lambda(v,\by\cup u;\theta)}{\lambda(v,\by;\theta)}\right\}.
\]
Once we have obtained the function $\widetilde \phi$, we obtain the
semi-optimal function $\phi$ by
$\phi(u,\by;\theta)=\widetilde\phi(u,\by;\theta)/\sqrt{\lambda(u,\by;\theta)}$.

\subsection{Numerical solution using Nystr{\"om} approximation}

The equation \eqref{eq:phiOpt2} is solved numerically using the
Nystr{\"om} approximation \citep{nystrom:1930}. We introduce a quadrature scheme with $m$
quadrature points \linebreak$u_1,\dots,u_m \in W$ and associated weights
$w_j$, $j=1,\ldots,m$, and approximate the operator $T_\by$ for any $\R^p$ valued function $g$ by
\begin{equation*}
T_\by g(u) \approx \sum_{j=1}^m  g(u_j)t(u,u_j,\by;\theta) w_j.
\end{equation*}
Introducing the quadrature approximation in \eqref{eq:phiOpt2}  and
multiplying each term by $\sqrt{w_i}$, we obtain $\sqrt{w_i}\widetilde
\phi(u_i,\by;\theta)$, $i=1,\ldots,m$, as solutions of the linear equations
\[
  \sqrt{w_i}\widetilde\phi(u_i,\by;\theta) + \sum_{j=1}^m
  \sqrt{w_iw_j}\widetilde t(u_i,u_j,\by;\theta) \sqrt{w_j}\widetilde\phi(u_j,\by;\theta) =
  \sqrt{w_i} \; \frac{\lambda^{(1)}(u_i,\by;\theta)}{\sqrt{\lambda(u_i,\by;\theta)}},
\]
for $i=1,\ldots,m$. These equations can be reformulated as the matrix equation
\begin{equation}\label{eq:matrixeq}
  \left\{\bI_{m}+\widetilde\bT(\by;\theta) \right\}  \sqrt{w} \widetilde\phi(\by;\theta) =  \ell(\by;\theta)
\end{equation}
where $\bI_m$ is the $m \times m$ identity matrix, $\widetilde{\bT} (\by;\theta)  =  \left\{ \sqrt{w_i w_j}
  \widetilde t(u_i,u_j,\by;\theta)\right\}_{ij}$, $i,j=1,\dots,m$, $\ell(\by;\theta)$ is the $m \times p$ matrix with rows $\sqrt{w_i} \lambda^{(1)}(u_i,\by;\theta)^\top/\sqrt{\lambda(u_i,\by;\theta)}$, $i=1,\dots,m$,
and $\sqrt{w}\widetilde \phi(\by;\theta)$ is the
$m\times p$ matrix with rows $\sqrt{w_i}\widetilde \phi(u_i,\by;\theta)^\top$.
The symmetric matrix $\widetilde\bT(\by;\theta)$ is sparse due to the finite
range property. Thus,
provided that $\bI_{m}+\widetilde\bT(\by;\theta)$ is positive definite, the
matrix equation can be solved with respect to $\sqrt{w}\widetilde \phi(\by;\theta)$ using sparse Cholesky factorization (see
\nocite{davis:06} Davis, 2006, and the \texttt R package \texttt{Matrix}).

Having solved \eqref{eq:matrixeq} with respect to
$\sqrt{w}\widetilde \phi(\by;\theta)$, and thus obtaining estimates of \linebreak$\sqrt{w_i} \widetilde
\phi(u_i,\by;\theta)$, we obtain estimates $\hat \phi(u_i,\by;\theta)$ of
$\phi(u_i,\by;\theta)$ via the
relation\linebreak $\hat\phi(u_i,\by;\theta) = \widetilde
\phi(u_i,\by;\theta)/\sqrt{w_i \lambda(u_i,\by;\theta)}$. Letting $\hat \phi(\by;\theta)$ be the $m
  \times p$ matrix with rows $\hat \phi(u_j,\by;\theta)^\top$ the
    Nystr{\"o}m approximation of  $\phi(u,\by;\theta)$ for any $u \in W$
is
\[ \hat \phi(u,\by;\theta) \approx
\frac{\lambda^{(1)}(u,\by;\theta)}{\lambda(u,\by;\theta)} -
\hat \phi(\by;\theta)^\top \{w_j t(u,u_j,\by;\theta)\}_{j=1}^m . \]
In particular, we obtain the approximations $\hat \phi(u,\bx \setminus u;\theta)$ of $\phi(u,\bx \setminus
u;\theta)$, $u \in \bx$, which are needed to evaluate the first term in
\eqref{eq:optimalef}. Finally, the
integral term in \eqref{eq:optimalef} and the empirical sensitivity
are approximated by
\[  \hat \phi(\bx;\theta)^\top \{w_j \lambda(u_j,\bx;\theta)\}_{j=1}^m \quad \text{ and }\quad
 \hat \phi(\bx;\theta)^\top w\lambda^{(1)}(\bx;\theta) \]
where $w\lambda^{(1)}(\bx;\theta) $ is the $m \times p$ matrix with rows $w_j \lambda^{(1)}(u_j,\bx;\theta)^\top$.

\subsection{Some computational considerations}
The matrix $\bI_{m}+\widetilde\bT(\by;\theta)$ is not guaranteed to be positive
definite. In case of purely repulsive point processes (Papangelou
conditional intensity always decreasing when neighbouring points are
added), all entries in $\widetilde\bT(\by;\theta)$ are positive and we did not
experience negative definite $\bI_{m}+\widetilde\bT(\by;\theta)$. However, with models
allowing for positive interaction, we occasionally experienced negative
definiteness in which case a solution for $\phi(\cdot,\by;\theta)$ cannot be
obtained. In such case we simply returned the pseudolikelihood
estimate.

In case of a quadrature scheme corresponding to a
subdivision of $W$ into square cells of sidelength $s$ the computational complexity of one Newton-Raphson update is
roughly of the order $(n+1) m (R/s)^2$. Thus, the semi-optimal approach is less feasible for data with a high number $n$ of points. 
Thus we would not recommend using the method for datasets with thousands of points while it is presently quite feasible for datasets with a few hundred points.

In case of e.g.\ the Strauss hard core process we may encounter $\lambda(u_j,\by;\theta)=0$. In this case, we use the conventions
$\lambda^{(1)}(u_j,\by;\theta)/\sqrt{\lambda(u_j,\by;\theta)}=0$
and $\lambda(u_j,\by \cup u_j;\theta)/\lambda(u_j,\by;\theta)=0$.

\end{document}